\documentclass[12pt]{article}
\usepackage{amsmath,amssymb}
\usepackage{mathrsfs}
\usepackage{amsfonts}
\usepackage{latexsym}
\usepackage{amsthm}
\usepackage{pictex}

\usepackage[T2A]{fontenc}     
\usepackage[cp1251]{inputenc} 
\usepackage{graphicx}

\newcommand{\er}[1]{{\rm(\ref{#1})}}
\def\lb{\label}
\theoremstyle{plain}
\newtheorem{theorem}{\bf Theorem}[section]
\newtheorem{lemma}[theorem]{\bf Lemma}
\theoremstyle{remark}

\begin{document}

\def\a{\alpha} \def\cA{{\cal A}} \def\bA{{\bf A}}  \def\mA{{\mathscr A}}
\def\b{\beta}  \def\cB{{\cal B}} \def\bB{{\bf B}}  \def\mB{{\mathscr B}}
\def\g{\gamma} \def\cC{{\cal C}} \def\bC{{\bf C}}  \def\mC{{\mathscr C}}
\def\G{\Gamma} \def\cD{{\cal D}} \def\bD{{\bf D}}  \def\mD{{\mathscr D}}
\def\d{\delta} \def\cE{{\cal E}} \def\bE{{\bf E}}  \def\mE{{\mathscr E}}
\def\D{\Delta} \def\cF{{\cal F}} \def\bF{{\bf F}}  \def\mF{{\mathscr F}}
\def\c{\chi}   \def\cG{{\cal G}} \def\bG{{\bf G}}  \def\mG{{\mathscr G}}
\def\z{\zeta}  \def\cH{{\cal H}} \def\bH{{\bf H}}  \def\mH{{\mathscr H}}
\def\e{\eta}   \def\cI{{\cal I}} \def\bI{{\bf I}}  \def\mI{{\mathscr I}}
\def\p{\psi}   \def\cJ{{\cal J}} \def\bJ{{\bf J}}  \def\mJ{{\mathscr J}}
\def\vT{\Theta}\def\cK{{\cal K}} \def\bK{{\bf K}}  \def\mK{{\mathscr K}}
\def\k{\kappa} \def\cL{{\cal L}} \def\bL{{\bf L}}  \def\mL{{\mathscr L}}
\def\l{\lambda}\def\cM{{\cal M}} \def\bM{{\bf M}}  \def\mM{{\mathscr M}}
\def\L{\Lambda}\def\cN{{\cal N}} \def\bN{{\bf N}}  \def\mN{{\mathscr N}}
\def\m{\mu}    \def\cO{{\cal O}} \def\bO{{\bf O}}  \def\mO{{\mathscr O}}
\def\n{\nu}    \def\cP{{\cal P}} \def\bP{{\bf P}}  \def\mP{{\mathscr P}}
\def\r{\rho}   \def\cQ{{\cal Q}} \def\bQ{{\bf Q}}  \def\mQ{{\mathscr Q}}
\def\s{\sigma} \def\cR{{\cal R}} \def\bR{{\bf R}}  \def\mR{{\mathscr R}}
\def\S{\Sigma} \def\cS{{\cal S}} \def\bS{{\bf S}}  \def\mS{{\mathscr S}}
\def\t{\tau}   \def\cT{{\cal T}} \def\bT{{\bf T}}  \def\mT{{\mathscr T}}
\def\f{\phi}   \def\cU{{\cal U}} \def\bU{{\bf U}}  \def\mU{{\mathscr U}}
\def\F{\Phi}   \def\cV{{\cal V}} \def\bV{{\bf V}}  \def\mV{{\mathscr V}}
\def\P{\Psi}   \def\cW{{\cal W}} \def\bW{{\bf W}}  \def\mW{{\mathscr W}}
\def\o{\omega} \def\cX{{\cal X}} \def\bX{{\bf X}}  \def\mX{{\mathscr X}}
\def\x{\xi}    \def\cY{{\cal Y}} \def\bY{{\bf Y}}  \def\mY{{\mathscr Y}}
\def\X{\Xi}    \def\cZ{{\cal Z}} \def\bZ{{\bf Z}}  \def\mZ{{\mathscr Z}}
\def\O{\Omega}
\def\ve{\varepsilon}
\def\vt{\vartheta}
\def\vp{\varphi}
\def\vk{\varkappa}

\def\mM{M}
\def\mB{B}
\def\mR{R}

\def\mA{{\mathscr A}}
\def\mB{{\mathscr B}}
\def\mC{{\mathscr C}}
\def\mD{{\mathscr D}}
\def\mE{{\mathscr E}}
\def\mF{{\mathscr F}}
\def\mG{{\mathscr G}}
\def\mH{{\mathscr H}}
\def\mI{{\mathscr I}}
\def\mJ{{\mathscr J}}
\def\mK{{\mathscr K}}
\def\mL{{\mathscr L}}
\def\mM{{\mathscr M}}
\def\mN{{\mathscr N}}
\def\mO{{\mathscr O}}
\def\mP{{\mathscr P}}
\def\mQ{{\mathscr Q}}
\def\mR{{\mathscr R}}
\def\mS{{\mathscr S}}
\def\mT{{\mathscr T}}
\def\mU{{\mathscr U}}
\def\mV{{\mathscr V}}
\def\mW{{\mathscr W}}
\def\mX{{\mathscr X}}
\def\mY{{\mathscr Y}}
\def\mZ{{\mathscr Z}}

\def\Z{{\Bbb Z}}
\def\R{{\Bbb R}}
\def\C{{\Bbb C}}
\def\T{{\Bbb T}}
\def\N{{\Bbb N}}
\def\S{{\Bbb S}}
\def\H{{\Bbb H}}
\def\J{{\Bbb J}}

\def\qqq{\qquad}
\def\qq{\quad}
\newcommand{\ma}{\begin{pmatrix}}
\newcommand{\am}{\end{pmatrix}}
\newcommand{\ca}{\begin{cases}}
\newcommand{\ac}{\end{cases}}
\let\ge\geqslant
\let\le\leqslant
\let\geq\geqslant
\let\leq\leqslant
\def\ma{\left(\begin{array}{cc}}
\def\am{\end{array}\right)}
\def\iint{\int\!\!\!\int}
\def\lt{\biggl}
\def\rt{\biggr}
\let\geq\geqslant
\let\leq\leqslant
\def\[{\begin{equation}}
\def\]{\end{equation}}
\def\wh{\widehat}
\def\wt{\widetilde}
\def\pa{\partial}
\def\sm{\setminus}
\def\es{\emptyset}
\def\no{\noindent}
\def\ol{\overline}
\def\iy{\infty}
\def\ev{\equiv}
\def\/{\over}
\def\ts{\times}
\def\os{\oplus}
\def\ss{\subset}
\def\h{\hat}
\def\Re{\mathop{\rm Re}\nolimits}
\def\Im{\mathop{\rm Im}\nolimits}
\def\supp{\mathop{\rm supp}\nolimits}
\def\sign{\mathop{\rm sign}\nolimits}
\def\Ran{\mathop{\rm Ran}\nolimits}
\def\Ker{\mathop{\rm Ker}\nolimits}
\def\Tr{\mathop{\rm Tr}\nolimits}
\def\const{\mathop{\rm const}\nolimits}
\def\dist{\mathop{\rm dist}\nolimits}
\def\diag{\mathop{\rm diag}\nolimits}
\def\Wr{\mathop{\rm Wr}\nolimits}
\def\BBox{\hspace{1mm}\vrule height6pt width5.5pt depth0pt \hspace{6pt}}

\def\Diag{\mathop{\rm Diag}\nolimits}

\def\Twelve{
\font\Tenmsa=msam10 scaled 1200 \font\Sevenmsa=msam7 scaled 1200
\font\Fivemsa=msam5 scaled 1200 \textfont\msbfam=\Tenmsb
\scriptfont\msbfam=\Sevenmsb \scriptscriptfont\msbfam=\Fivemsb

\font\Teneufm=eufm10 scaled 1200 \font\Seveneufm=eufm7 scaled 1200
\font\Fiveeufm=eufm5 scaled 1200
\textfont\eufmfam=\Teneufm \scriptfont\eufmfam=\Seveneufm
\scriptscriptfont\eufmfam=\Fiveeufm}

\def\Ten{
\textfont\msafam=\tenmsa \scriptfont\msafam=\sevenmsa
\scriptscriptfont\msafam=\fivemsa

\textfont\msbfam=\tenmsb \scriptfont\msbfam=\sevenmsb
\scriptscriptfont\msbfam=\fivemsb

\textfont\eufmfam=\teneufm \scriptfont\eufmfam=\seveneufm
\scriptscriptfont\eufmfam=\fiveeufm}

\title {Borg type uniqueness Theorems for periodic Jacobi operators
with matrix valued coefficients}

\author{
 Evgeny Korotyaev
\begin{footnote}
{ Institut f\"ur  Mathematik,  Humboldt Universit\"at zu Berlin,
Rudower Chaussee 25, 12489, Berlin, Germany, e-mail:
evgeny@math.hu-berlin.de}
\end{footnote}
 \and Anton Kutsenko
\begin{footnote}
{ Department of
 Mathematics of Sankt-Petersburg State University, Russia e-mail: kucenkoa@rambler.u}
\end{footnote}
}

\maketitle

\begin{abstract}
\no
We give a simple proof of
Borg type uniqueness Theorems for periodic Jacobi operators
with matrix valued coefficients.
\end{abstract}


\section {Introduction}
\setcounter{equation}{0} Consider a self-adjoint matrix-valued
Jacobi operator $J$ acting on $\ell^2(\Z)^m$ and given by
\[
 \lb{000}
 (J y)_n=a_n y_{n+1}+b_ny_n+a_{n-1} y_{n-1},\qq n\in\Z,\qq
 y_n\in \C^m, y=(y_n)_{n\in \Z}\in \ell^2(\Z)^m,m\ge 1,
\]
 where $a_n>0,b_n=b_n^*$, $n\in\Z$ are  p-periodic sequences of the complex
 $m\ts m$ matrices. It is well known that the spectrum $\s(J)$ of $J$ is absolutely
continuous  and consists of non-degenerated intervals $[\l^+_{n-1},\l^-_n], \l^+_{n-1}<\l^-_n \le  \l^+_{n}, n=1,...,N<\iy$. 
These intervals are separated by the gaps $\g_n=(\l_n^-,\l_n^+),
n=1,..,N-1$ with the length $>0$. Introduce the fundamental $m\ts m$
matrix-valued solutions $\vp=(\vp_n(z))_{n\in\Z},
\vt=(\vt_n(z))_{n\in\Z}$ of the equation
\[
 \lb{001}
 a_ny_{n+1}+b_ny_n+a_{n-1}^* y_{n-1}=z y_n,\qqq
\vp_{0}\ev \vt_1\ev 0,\ \vp_1\ev \vt_{0}\ev I_m,\qq
 (z,n)\in\C\ts\Z,
\]
where $I_m$ is the identity $m\ts m$ matrix. We define the
monodromy $2m\ts2m$ matrix $\mM_p$  by
\[
 \lb{002}
 \mM_p(z)=\ma \vt_{p}(z) & \vp_{p}(z) \\
           \vt_{p+1}(z) & \vp_{p+1}(z)
          \am,\qq z\in \C.
\]
Let $\t_1(z),..,\t_{2m}(z)$ be eigenvalues of $\mM_p(z)$. Recall
that $\s(J)=\bigcup_{j=1}^{2m}\{z\in \C:\ |\t_j(z)|=1\}$, see
\cite{KKu}. Let $J^0$ be the unperturbed Jacobi matrix  with
$a_n^0=I_m$, $b_n^0=0$. Note, that $\s(J^0)=\s_{ac}(J^0)=[-2,2]$.
Let $\N_s=\{1,..,s\}$ and let $c=\det\prod_{n=1}^pa_n$. We formulate
first Borg type uniqueness theorem

\begin{theorem}
\lb{t030}
Let $c=1$ and let $
 \s(J)=[\l_0^+,\l_N^-]=\{each \ |\t_j(z)|=1,\ (z,j)\in
 [\l_0^+,\l_N^-]\ts \N_{2m}\}$, where $\l_0^+=-\l_N^-$.
Then $J=J^0$.
\end{theorem}

This theorem was proved in \cite{CGR}. We give a simple proof based
on the trace formula and the properties of the Chebyshev polynomials
proved in Lemma \ref{t010}. There is an enormous literature on
inverse spectral problems for scalar (i.e., $m=1$) periodic Jacobi
matrices (see \cite{BGGK}, \cite{K},\cite{KKu1}, \cite{KKu2},
\cite{vM}, \cite{T}, book \cite{T} and references therein), but very
little for matrix-valued periodic Jacobi operators (see \cite{CGR}
and references therein). Note that the complete solution of  inverse
problem for finite matrix-valued Jacobi operators was obtained
recently \cite{BCK}.

The operator $J$ is unitarely equivalent to the operator
$\mJ=\int_{[0,2\pi)}^{\os}K_p(e^{ix}){dx\/2\pi}$ acting in
$\int_{[0,2\pi)}^{\os}\mH{dx\/2\pi}$, where $\mH=\C^p$ and  $pm\ts
mp$ matrix $K_p(\t)$ is given by
\[
 \lb{t001}
 K_p(\t)=\left(\begin{array}{ccccc} b_1 & a_1 & 0 & ... & \t^{-1} a_p \\
                                  a_1 & b_2 & a_2 & ... & 0 \\
                                  0 & a_2 & b_3 & ... & 0 \\
                                  ... & ... & ... & ... & ... \\
                                  \t a_p & 0 & ... & a_{p-1} & b_p
                                  \end{array}\right),\ \ \
                                  \t\in \S^1=\{\t:|\t|=1\},
\]
see [KKu]. Let $(\l_{n}(\t))_{n=1}^{pm}$ be eigenvalues of
$K_p(\t)$.   We formulate our main result.

\begin{theorem} \lb{t031}
Let $a_n>0,b_n=b_n^*$ be real $m\ts m$ matrices for all $n\in\Z$.

\no i) Let $c=1$. Then $\sum_{n=1}^{mp}\l_n^2(\t)=2pm$ for some
$\t\in \S^1$ iff $J=J^0$.

\no ii)  Let $c=1$ and let $\vk_1\in\R$.  Then eigenvalues
$\l_{s}(e^{i\vk_1})=2\cos\frac{\vk_1+2\pi(s-1)}p$ for all
$s\in\N_{mp}$ iff $J=J^0$.

\no iii) Let $\vk_1,\vk_2\in\R$ and $\cos\vk_1\not=\cos\vk_2$. Let
eigenvalues $\l_{s}(e^{i\vk_1})=2\cos\frac{\vk_1+2\pi(s-1)}p$ for
all $s\in\N_{mp}$ and $\l_n(e^{i\vk_2})=2\cos {\vk_2+2\pi n_1\/p},
n\in\N_m$, for some $n_1\in\N_{p}$, and if $e^{i2\vk_2}=1$, then
additional eigenvalues $\l_{n+m}(e^{i\vk_2})=\l_1(e^{i\vk_2})\ne \pm
2$, $n\in\N_m$. It is possible iff $J=J^0$.

\end{theorem}


{\no\bf Remark.} 1) Note that the condition  $\cos\frac{\vk_2+2\pi
n_1}p\not=\pm1$ in iii) is associated with the
 unperturbed operator $J^0$, where the endpoints of the spectrum
 $\s(J^0)=[-2,2]$ have the multiplicity $m$,
 as the zeros of the determinant $D_p(z,\pm 1)$. Each  point from $(-2,2)$
 has the multiplicity $2m$ as  the zero of the determinant $D_p(z,\t), \t\ne \pm 1$.

\no 2) Consider the case $m=2$ and $\vk_1=0$, $\vk_2=\pi$. Let
periodic eigenvalues $\l_{s}(1)=2\cos\frac{\vk_1+2\pi(s-1)}p$ for
all $s\in\N_{2p}$ and let anti-periodic eigenvalues
$\l_s(-1)=2\cos{{\pi(2n+1)}\/p}\ne\pm2, s\in\N_4$ for some
$n\in\N_{p}$. Then we deduce that $J=J^0$.


\no 3) Consider a self-adjoint matrix-valued Jacobi operator $\wt J$
acting on $\ell^2(\Z)^m$ and given by
$$
 (\wt J y)_n=\wt a_n y_{n+1}+\wt b_ny_n+\wt a_{n-1}^* y_{n-1},\qq n\in\Z,\qq
 y_n\in \C^m, y=(y_n)_{n\in \Z}\in \ell^2(\Z)^m,m\ge 1,
$$
 where $\wt a_n$ ($\det \wt a_n\ne 0$), $\wt b_n=\wt b_n^*$, $n\in\Z$ are
 p-periodic sequences of the complex
 $m\ts m$ matrices. Then the operator $\wt J=UJ U^*$
 where $J$ is given by \er{001} and $U=\diag_{n\in\Z}u_n$
and the unitary matrices $u_n$ have the forms $u_0=I_m$;
$u_{n+1}=c_n^*u_n$, $n\ge0$; $u_{n-1}=c_{n-1}u_{n+1}$, $n\le0$. Here
$\wt a_n=c_n a_n, a_n>0$ and $c_n$ are the unitary matrices.

\section{Proof}
\setcounter{equation}{0}

We need the following results from \cite{KKu}.

\begin{lemma}\lb{7100}\it
\no The following identities and estimate are fulfilled:
\[
 \lb{2010}
 \sum_{n=1}^{mp}\l_{n,p}(\t)=\sum_{n=1}^p\Tr b_n,\qqq
 \sum_{n=1}^{mp}\l_{n,p}^2(\t)=\sum_{n=1}^p\Tr (b_n^2+2a_n^2),
 \qq \  all \qq \t\in\C\sm\{0\},
\]
\[
 \lb{7110}
 \sum_{n=1}^{pm} \l_{n,p}^2(\t)\ge2pm c^{2\/{pm}},\qqq \t\in \S^1,
\]
where the identity \er{7110} holds true iff $J=c^{1\/pm}J^0$, and
here $c=\det \prod_1^p a_n$.
\end{lemma}

Recall the following identities for $D_p(z,\t)=\det (\mM_p(z)-\t
I_{2m})$ from \cite{KKu}
\begin{multline}
 \lb{8000}
 D_p(z,\t)=\prod_{j=1}^{2m}(\t_{j}(z)-\t)=
 c(-\t)^m\prod_{n=1}^{mp}(z-\l_{n}(\t))=c(-\t)^m\det(z-K_p(\t)), \
\end{multline}
where $z,\t\in\C,\qq \t\ne 0$.

{{\no\bf Proof of Theorem \ref{t030}}}  {\bf Consider the first
case} $\sum_{n=1}^p\Tr b_n=0$. Consider the $p$-periodic operator
$J$ as $pk$-periodic operator for some $(k,\t)\in \N\ts \S^1$. Let
the eigenvalues of
 $K_{pk}(\t)$ be given by
\[
 \lb{o000}
 \l_-(\t)\ev\l_{1,pk}(\t)\le\l_{2,pk}(\t)\le..\le\l_+(\t)\ev\l_{pmk,pk}(\t).
\]
The eigenvalues $\l_{\pm}(\t)$ belong to $\s(J)$, then we have that
$[\l_-(\t),\l_+(\t)]\ss\s(J)$. Using $|\t_j(z)|=1$ for all
$(z,j)\in[\l_-(\t),\l_+(\t)]\ts\N_{2m}$, and \er{8000}, we obtain
$$
 |D_{pk}(z,\t)|=\left|\prod_{j=1}^{2m}(\t-\t_j^k(z))
 \right|\le\prod_{j=1}^{2m}(|\t|+|\t_j^k(z)|)\le2^m,\ \ all\ \
 z\in[\l_-(\t),\l_+(\t)].
$$
We have $\det\prod_{n=1}^{pk}a_n=c^k=1$, since a sequence $a_n$ is
$p$-periodic. Then \er{8000} gives
\[
 \lb{u010}
 |D_{pk}(z,\t)|=\left|\prod_{n=1}^{pmk}(z-\l_{n,pk}(\t))\right|\le2^m,\
 \ all \  \ z\in[\l_-(\t),\l_+(\t)].
\]
A sequence $b_n$ is  $p$-periodic, then
\[
 \lb{u011}
 \sum_{n=1}^{pmk}\l_{n,pk}(\t)=\sum_{n=1}^{pk}\Tr b_n=k\sum_{n=1}^p\Tr b_n=0.
\]
The relations \er{o000}, \er{u010}, \er{u011} and Lemma \ref{t010}
give that
\[
 \lb{t020}
 (\l_{n,pk})_{n=1}^{mpk}\in\cP_{mpk}(2^m),\qqq
 \sum_{n=1}^{pmk}\l_{n,pk}^2(\t)\le 2pmkC_k, \qqq
 C_k=2^{\frac{2(m-1)}{pmk}}.
\]
Then using \er{2010} and  $a_n=a_{n+p}, b_n=b_{n+p}$ for all $n\in
\Z$, we obtain
$$
 \sum_{n=1}^{pmk}\l_{n,pk}^2(\t)=\sum_{n=1}^{pk}\Tr(b_n^2+a_n^2)=
 k\sum_{n=1}^{p}\Tr(b_n^2+a_n^2)=k\sum_{n=1}^{pm}\l^2_{n,p}\le
 k(2pm)C_k,
$$
which yields $\sum_{n=1}^p\l^2_{n,p}\le2pm$, since  we  take
$k\to+\iy$. Thus Lemma \ref{7100} implies $J=J^0$.

{\bf Consider the second case} $\b=\sum_{n=1}^p\Tr b_n$. Then a new
operator $J^1=(J-\b)$ satisfies $\sum_{n=1}^p\Tr (b_n-\b I_m)=0$.
Thus due to the first case we deduce that
$$
\s(J)=\s(J^1)=[-2,2]=[-x-\b,x-\b],\qq a_n=I_m, \qq b_n-\b
I_m=b_n^1=0
$$
 for $x=\l_n^-$, which yields
$x=2$, $\b=0$ and the theorem has been proved.
\BBox

\begin{lemma} \lb{o001}
i) Let $\l$ be an eigenvalue of $K_p(\t)$ and have multiplicity $m$
for some $\t\in\S^{1}\sm\{-1,1\}$. Then the multiplaers have the
form $\t_j(\l)=\t, \t_{j+m}(\l)=\t^{-1}$ for all $j\in\N_{m}$.

\no ii) Let $\l$ be an eigenvalue of $K_p(\t)$ and have multiplicity
$2m$ for some $\t\in\{-1,1\}$. Then each $\t_j(\l)=\t$,
$j\in\N_{2m}$.
\end{lemma}
\no {\bf Proof.} i) The matrix $K_p(\t)$ is self-adjoint. Let
$K_p(\t)f^k=\l f^k$ for some  orthogonal eigenvectors
$f^k=(f^k_n)_{n=1}^{mp}$, $k\in\N_m$. If $f^k_0=\t^{-1}f^k_p$, then
the definition of the matrix $\mM_p$ gives $
\mM_p(\l)(f^k_0,f^k_1)^\top=\t(f^k_0, f^k_1)^\top$. Note that $\wt
f^k=(f^k_0,f^k_1)^\top$ define other components of vector $f^k$,
since $L(\t)$ has special form, see \er{t001}. Then the vectors $\wt
f^k$, $k\in\N_m$ are linearly independent vectors, since $f^k$ are
linearly independent vectors. Then $\t$ has multiplicity at least
$m$. The matrix $\mM_p$ is symplectic, then $\t^{-1}$ is eigenvalue
of $\mM_p$ and $\t^{-1}$ has multiplicity at least $m$. We obtain
first statement, since $\mM_p$ is $2m\ts 2m$ matrix. The proof of
ii) is similar. \BBox

{{\no\bf Proof of Theorem \ref{t031}.}} The statement i) follows
from Lemma \ref{7100}.

ii) {\bf Sufficiency}. 
Recall that if $m=1$, then $\det (M_p^0-\t I_{2})=
 \t(\t+\t^{-1}-2\mT_p(z/2))$, where $\mT_p(z)=\cos(p\arccos z)$ is the Chebyshev
polynomial. Moreover, zeros of the polynomial $\mT_p(z)-\cos \vk_1$
are given by $\l_{s}(e^{i\vk_1})=2\cos\frac{\vk_1+2\pi(s-1)}p,
s\in\N_{p}$.

Thus if $J=J^0, m\ge 2$, then corresponding monodromy operators
$M_p$ satisfies $ \det (M_p-\t I_{2m})=
 \t^{m}\prod_{j=1}^m(\t+\t^{-1}-2\mT_p(z/2)),$
 which yields at $\t=e^{i\vk_1}$
\[
\lb{imk}
det (M_p-e^{i\vk_1} I_{2m})=
 \t=e^{im\vk_1}2^m(\cos \vk_1-\mT_p(z/2))^m=e^{im\vk_1}2^m,
\]
This implies
$\l_{s}(e^{i\vk_1})=2\cos\frac{\vk_1+2\pi(s-1)}p, s\in\N_{mp}$.

{\bf Necessity.} Let
$\l_{s}(e^{i\vk_1})=2\cos\frac{\vk_1+2\pi(s-1)}p, s\in\N_{mp}$. Then
the direct calculation implies
$\sum_{n=1}^{mp}\l_n^2(e^{i\vk_1})=2pm$, and the statement i) gives
$J=J^0$.

iii) {\bf The sufficiency} is proved similar to the case ii).

{\bf Necessity}. We consider only the case $e^{2i\vk_2}\not=1$, the
proof of the case $e^{2i\vk_2}\not\ne 1$ is similar. Let
$\l_{s}(e^{i\vk_1})=2\cos\frac{\vk_1+2\pi(s-1)}p, s\in\N_{mp}$.
Using \er{8000}, \er{imk}, we obtain
\[
 \lb{a0002}
 \prod_{j=1}^{2m}(e^{i\vk_1}-\t_j(z))=c(-1)^m{e^{im\vk_1}}
 \prod_{n=1}^{mp}(z-\l_{n}(e^{i\vk_1}))=
 c(-2)^m{e^{im\vk_1}}(\mT_p(z/2)-\cos\vk_1)^m,
\]
where $c=\det\prod_{n=1}^pa_n $. The eigenvalue $\l_1(e^{i\vk_2})$
has the multiplicity $m$, then using Lemma \ref{o001}  and
substituting $z=\l_1(e^{i\vk_2})=2\cos\frac{\vk_2+2\pi n_1}p$ and
two multipliers $e^{\pm i\vk_2}$  (given by Lemma \ref{o001}) into
\er{a0002}, we obtain
$$
 (e^{i\vk_1}-e^{i\vk_2})^m(e^{i\vk_1}-e^{-i\vk_2})^m=
 c2^m{e^{im\vk_1}}(\cos\vk_2-\cos\vk_1)^m,
$$
which yields $c=1$.  Then the statement ii) gives $J=J^0$.
 \BBox

\begin{lemma}
\lb{t010} For any $r\ge0$, $s\ge2$ the following identity holds true
$$
\sup_{x\in\cP_s(r)}\sum_1^s x_n^2=2s\left(\frac r2\right)^{\frac2s},
$$
$$
 \cP_s(r)=\left\{(x_n)_{1}^s:\ x_1\le..\le x_s,\ \sum_{n=1}^s x_n=0,
 \ \left|\prod_{n=1}^s(z-x_n)\right|\le r,\ \ all \
 z\in[x_1,x_s]\right\}\ss\R^s.
$$
\end{lemma}
{\no\bf Proof.} In the proof we use arguments from \cite{Ku}. Let
$\|x\|^2=\sum x_n^2, x=(x_n)_1^s\in \R^s$. The set $\cP_s(r)$ is
compact, then
$ \sup_{x\in\cP_s(r)}\|x\|^2=\|x^0\|^2$ for some $
x^0=(x_n^0)_{1}^s\in\cP_s(r)$. Introduce a polynomial
$p_0(z)=\prod_{n=1}^s(z-x^0_n)$. Let $x_0^0=-\iy$, $x_{s+1}^0=+\iy$.
The polynomial $p_0'$ has only real zeros $\wt x_n, n\in \N_{s-1}$.
We will show that each $|p_0(\wt x_n)|=r, n\in \N_{s-1}$. Assume
that there exist $1\le n_1<n_2\le s$ such that
\[
 \lb{u002}
 x^0_{n_1-1}<x^0_{n_1}\le x^0_{n_2}<x^0_{n_2+1},\ \
 \max_{z\in[x_{n_1},x_{n_2}]}|p_0(z)|<r.
\]
Introduce a polynomial $p_{\ve}(z)=\prod_{n=1}^s(z-x_n^{\ve})$,
where
 a vector $x^{\ve}=(x^{\ve}_n)_{n=1}^s\in\R^s$ is given by
\[
 \lb{u003}
 x^{\ve}_n=x^0_n,\ n\not=n_1,\ n\not=n_2;\ \ x^{\ve}_{n_1}=x^0_{n_1}-\ve,\
 x^{\ve}_{n_2}=x^0_{n_2}+\ve,  \ \ \ve>0.
\]
Using $x^0\in\cP_s(r)$ and \er{u002}, we obtain
\[
 \lb{u004}
 x_1^{\ve}\le..\le x_s^{\ve},\ \ \sum_{n=1}^s
 x_n^{\ve}=\sum_{n=1}^sx^0_n=0,\ \
 \max_{z\in[x^{\ve}_{n_1},x^{\ve}_{n_2}]}|p_{\ve}(z)|\le r,
\]
for sufficiently small $\ve>0$. We rewrite $p_{\ve}$ in the form
\[
 \lb{u005}
 p_{\ve}(z)=p_0(z)g_{\ve}(z),\ \ g_{\ve}(z)=
 \frac{(z-x^{\ve}_{n_1})(z-x^{\ve}_{n_2})}
 {(z-x^0_{n_1})(z-x^0_{n_2})}=
 \frac{(z-x^{0}_{n_1}+\ve)(z-x^{0}_{n_2}-\ve)}
 {(z-x^0_{n_1})(z-x^0_{n_2})}.
\]
Due to $x^{0}_{n_1}\le x^{0}_{n_2}$ and $\ve>0$ we deduce that
$|g_{\ve}(z)|\le1$, $z\in\R\sm[x^{\ve}_{n_1},x^{\ve}_{n_2}]$. Then
\er{u005} yields $|p_{\ve}(z)|\le |p_0(z)|$,
$z\in\R\sm[x^{\ve}_{n_1},x^{\ve}_{n_2}]$ and thus  \er{u004} gives
$|p_{\ve}(z)|\le r$, $z\in[x^{\ve}_{1},x^{\ve}_{s}]$, and
$x^{\ve}\in\cP_s(r)$ for all sufficiently small $\ve>0$. There is an
estimate
\[
 \lb{u006}
 \|x^{\ve}\|^2=\sum_{n=1}^s(x^{\ve}_n)^2=
 \sum_{n=1}^s(x^{0}_n)^2+4\ve(x^0_{n_2}-x^0_{n_1})+2\ve^2>\|x^0\|^2,
\]
since $x^0_{n_1}\le x^0_{n_2}$ and $\ve>0$. But \er{u006} and the
condition $x^{\ve}\in\cP_s(r)$ contradict the identiy $
\sup_{x\in\cP_s(r)}\|x\|^2=\|x^0\|^2$. Then the assumption \er{u002}
is not true and each $|p_0(\wt x_n)|=r, n\in \N_{s-1}$. Note that
only polynomials  $r\mT_s(\a z+\b)$, $\a\in\R\sm\{0\}$, $\b\in\R$,
(here $\mT_s$ are the Chebyshev polynomials, i.e. $\mT_s(\cos
z)=\cos sz$) have  this property, then $p_0(z)=r\mT_s(\a z+\b)$ for
some $\a\in\R\sm\{0\}$, $\b\in\R$. We take $\a,\b$ such that
$p_0(z)=\prod_{n=1}^s(z-x_n^0)$, $\sum_{n=1}^s x_n^0=0$ and then
$p_0(z)=r\mT_s(z{r^{-\frac1s}2^{\frac{1-s}s}})$, since (see
\cite{AS})
$$
 \mT_s(z)=\frac12\sum_{k=0}^{[{s\/2}]}(-1)^k\frac{s}{s-k}C^{s-k}_k(2z)^{s-2k}=
 2^{s-1}z^s-2^{s-3}sz^{s-2}+o(z^{s-2})\qq as \
 z\to\iy,\
$$
where $C^n_m=\frac{n!}{m!(n-m)!}$. Then using Viette formulas we get
$$
 \mT_s(z)=2^{s-1}\prod_{n=1}^s(z-z_n)=2^{s-1}\left(z^s-\x z^{s-1}+
 \frac12\left(\x^2-\e\right)z^{s-2}+o(z^{s-2})\right)\qq as \
 z\to\iy,\
$$
where $z_n$, $n\in\N_s$ are zeroes of $\mT_s$ and $\x=\sum_1^sz_n$,
$\e=\sum_1^sz_n^2$. Then we get $\x=0$ and $\e=\frac s2$. This gives
$$
 \sup_{x\in\cP_s(r)}\|x\|^2=\|x^0\|^2=r^{\frac2s}2^{\frac{2s-2}s}\sum_1^s z_n^2=
 r^{\frac2s}2^{\frac{2s-2}s}\e=2s\left(\frac
 r2\right)^{\frac2s} \BBox
$$

 \no {\bf Acknowledgments.}
Evgeny Korotyaev was partly supported by DFG project BR691/23-1. The
some part of this paper was written at the Math. Institute of
Humboldt Univ., Berlin; Anton Kutsenko is grateful to the Institute
for the hospitality.

\end{document}